\documentclass[preprint,11pt]{elsart}

\usepackage{graphicx}
\usepackage{amssymb}
\usepackage{amssymb,amsfonts,amsmath}
\usepackage{epsfig}
\usepackage{graphicx,wrapfig}
\usepackage{color}
\newcommand{\RR}{{\mathbb{R}}}
\newcommand{\CC}{{\mathbb{C}}}

\newtheorem{lemma}{Lemma}% [section]
\newtheorem{theorem}{Theorem}
\definecolor{myred}{cmyk}{0.000000,1.000000,1.000000,0.0}
\definecolor{myblue}{cmyk}{1.000000,0.750000,0.000000,0.1}

\begin{document}

\begin{frontmatter}

\title{Revisiting the stability of computing the \\ roots of a quadratic polynomial}

\author{Nicola Mastronardi}
\address{Istituto per le Applicazioni del Calcolo ``M. Picone'', sede di Bari,
 Consiglio Nazionale delle Ricerche,
 Via G. Amendola, 122/D, I-70126 Bari, Italy.}
\ead{n.mastronardi@ba.iac.cnr.it}

\author{Paul Van Dooren} 
\ead{paul.vandooren@uclouvain.be}
\address{Catholic University of Louvain,
Department of Mathematical Engineering,
Avenue Georges Lemaitre 4,
B-1348 Louvain-la-Neuve, Belgium}

\begin{abstract}
We show in this paper that the roots $x_1$ and $x_2$ of a scalar quadratic polynomial $ax^2+bx+c=0$ with real or complex  coefficients $a$, $b$ $c$ can be computed in a element-wise mixed stable manner, measured in a relative sense. We also show that this is a stronger property than norm-wise backward stability, but weaker than element-wise backward stability. We finally show that there does not exist any method that can compute the roots in an element-wise backward stable sense, which is also illustrated by some numerical experiments.
\end{abstract}


\end{frontmatter}

\section{Introduction}

In this paper we consider the very simple problem of computing the two roots of a quadratic polynomial 
\begin{equation}\label{p} p(x) := ax^2+bx+c \end{equation}
where the coefficients $a,b,c$ are either in $\RR$ or in $\CC$ and where $a\neq 0$ in order the equation to have indeed two roots. This is a very classical problem for which the solution is well known, namely
$$ x_{1,2} = \frac{-b\pm \sqrt{b^2-4ac}}{2a}. $$
But the straightforward implementation of the above formula is quite often numerically unstable for special choices of the coefficients $a,b,c$. One would like, on the other hand, to have a computational scheme that produces the computed roots $\hat x_1$ and $\hat x_2$ which correspond to an element-wise backward stable error, i.e. the relative backward errors are of the order of the unit roundoff $u$ for each individual coefficient $a$, $b$ and $c$. In fact, we can assume that $a$ is not perturbed in this process. We will call this {\bf Element-wise Backward Stability (\bf EBS}) :
\begin{center}
$ a(x-\hat x_1)(x-\hat x_2)= ax^2+\hat b x + \hat c $ \\ \medskip
$ |b-\hat b| \leq \Delta |b|, \quad  |c-\hat c| \leq \Delta |c|, \quad \Delta = {\mathcal O}(u).$
\end{center}
We will see that this can not be proven in the general case, but instead, we can obtain the slightly weaker result of {\bf Element-wise Mixed Stability (\bf EMS}), which implies that the computed roots $\hat x_1$ and $\hat x_2$ satisfy
\begin{center} 
$ a(x-\tilde x_1)(x-\tilde x_2)= ax^2+\hat b x + \hat c $ \\ \medskip
$ |\hat x_1-\tilde x_1| \leq \Delta |\tilde x_1|, \;  |\hat x_2-\tilde x_2| \leq \Delta |\tilde x_2|,$  \\ \medskip
$ |b-\hat b| \leq \Delta |b|, \;  |c-\hat c| \leq \Delta |c|, \; \Delta = {\mathcal O}(u),$
\end{center}
which means that the computed roots are close to roots of a nearby polynomial, all in a relative element-wise sense.

This last property is also shown to be stronger than the so-called {\bf Norm-wise Backward Stability (NBS}) which only imposes that the vector of perturbed coefficients is close to the original vector in a relative norm sense~:
\begin{center}
$ a(x-\hat x_1)(x-\hat x_2)= ax^2+\hat b x + \hat c $ \\ \medskip
$ \left\| \left[ \begin{array}{ccc} a & b & c  \end{array}\right] -\left[ \begin{array}{ccc}  a &  \hat b & \hat c \end{array}\right ]\right\| \leq   \Delta \left\| \left[ \begin{array}{ccc}  a & \hat b & \hat c \end{array}\right]\right\| , \quad \Delta = {\mathcal O}(u).$
\end{center}

This problem was studied already by several authors, but we could not find any conclusive answer to the EBS of 
any of the proposed algorithms. 

In this paper, we will first consider the case of real coefficients since it is more commonly occurring and the results are slightly stronger. We then show how it extends to the case of complex coefficients. We end with a section on numerical experiments where we also show that there does {\bf not} exist a method that is {\bf EBS} for all quadratic polynomials.

\section{Real coefficients} 

Before handling the general case where all three coefficients are nonzero, we point out that when $b$ and/or $c$ are zero the proof of {\bf EBS} is rather simple. 

\subsection{A zero coefficient}
{\bf Case $c=0$}. \\
If $c=0$, then the roots can be computed as follows
$$ x_1:=-b/a, \quad x_2=0$$
which is element-wise backward stable since under the IEEE floating point standard, we have that the {\em computed} roots satisfy
$$ \hat x_1= -fl(b/a)=-b(1+\delta)/a=-\hat b/a, \quad \hat x_2 = 0, \quad |\delta| \leq u,$$
where $u$ is the unit round-off of the IEEE floating point standard (see \cite{Higham}).
The backward error then indeed satisfies the relative element-wise bounds
$$ |\hat b - b| \leq u |b|, \quad |\hat c - c| =0|c|. $$
{\bf Case $b=0$}. \\
If $b=0$ then the roots can be computed as follows
$$ \quad x_1 = \sqrt{-c/a}, \quad x_2 := -x_1,$$
which is also element-wise backward stable since under the IEEE floating point standard, we have that the  {\em computed} roots satisfy the element-wise bounds
$$ \hat x_1=fl(\sqrt{fl(-c/a)})= \sqrt{-c(1+\eta)/a}, \quad \hat x_2=-\hat x_1, \quad |\eta| \leq \gamma_3:=\frac{3u}{1-3u}. $$
Notice that if sign$(c)$=sign$(a)$, the roots are purely imaginary.
The backward error for this computation satisfies the relative element-wise bounds
$$ |\hat b - b| \leq 0|b|, \quad |\hat c - c| \leq \gamma_3 |c|. $$

\subsection{Preliminary scaling}

We can thus assume now that all coefficients are nonzero. We start by reducing the problem to a simpler ``standardized" form in order to simplify the computational steps.

{\bf Scaling the polynomial $p(x)$}\\  
We scale the polynomial coefficients so that it is monic : $b_1:=b/a, c_1:=c/a$, 
which can be performed in a backward and forward stable way since we assumed $a\neq 0$. According to the IEEE floating point standard we have that the computed values $\hat b_1=fl(b_1)$ and $\hat c_1=fl(c_1)$ satisfy
the relative element-wise bounds
$$ |b_1-\hat b_1|\leq u|b_1| , \; |c_1-\hat c_1|\leq u|c_1|. $$
This implies we can as well consider the monic polynomial $$p_1(x):=p(x)/a=x^2+b_1x+c_1.$$

\medskip

{\bf Scaling the variable $x$} \\ 
We transform the variable $x$ to $y:=-x/\alpha$ where $|\alpha|:=\sqrt{|c_1|}$ and $\mathrm{sign}(\alpha)=\mathrm{sign}(b_1)$, and consider the polynomial $p_1(-\alpha y)/\alpha^2$ which is now monic in $y$
\begin{equation}\label{y} q(y):=y^2 -2 \beta y +e =0, \end{equation}
and where $\beta \in \RR_+$ and $e=\pm 1$. The formulas to compute $\alpha$, $\beta$ and $e$ are 
$$ \alpha:=\mathrm{sign}(b_1)\sqrt{|c_1|}, \; \beta := |b_1|/(2\sqrt{|c_1|}), \; e:=\mathrm{sign}(c_1)\cdot 1.$$
Since the sign function is exact under relative perturbations, $e$ is computed exactly. It then follows
that $\alpha$ and $\beta$ can be performed in a backward and forward stable way~: the computed values $\hat \alpha=fl(\alpha)$ and $\hat \beta=fl(\beta)$ satisfy the relative element-wise bounds
$$ |\alpha-\hat \alpha|\leq u|\alpha| , \; |\hat \beta-\beta|\leq 2u|\beta|, $$
and $e$ is computed exactly.
This implies we can as well consider the polynomial $g(y)=y^2-2\beta y+e$. 
We recapitulate this in a formal lemma. 

\begin{lemma} The transformations 
$$ [\alpha,\beta]=g_a[b,c] \quad \mathrm{and} \quad [b,c]=g_a^{-1}[\alpha,\beta]$$ 
between the polynomial $p(x)=ax^2+bx+c, \; a\neq 0$ and the monic polynomial $p(-\alpha y)/(a\alpha^2)=q(y)=y^2-2\beta y +e$ defined by the forward and backward relations 
$$ \alpha:=\mathrm{sign}(b/a)\sqrt{|c/a|}, \quad \beta := |b/a|/(2\sqrt{|c/a|}),
$$
and
$$ b=-2a\beta\alpha, \quad c=ae\alpha^2,
$$
where $a$ and $e=\mathrm{sign}(c/a)\cdot 1$ are not perturbed, are both element-wise well-conditioned maps.
\end{lemma}
{\bf Proof}. 
If we define the perturbations for the forward map as $$[\alpha(1+\delta_\alpha),\beta(1+\delta_\beta)]=g_a[b(1+\delta_b),c(1+\delta_c)],$$ then the above discussion says that the relative perturbations $\delta_\alpha,\delta_\beta$ on the result are ${\cal O}(u)$ if the relative perturbations on the data $\delta_b,\delta_c$ are ${\cal O}(u)$. 
The same reasoning can be applied to the perturbation of the backward map  
$$[b(1+\delta_b),c(1+\delta_c)]=g_a^{-1}[\alpha(1+\delta_\alpha),\beta(1+\delta_\beta)],$$
which says now that 
$\delta_b,\delta_c={\cal O}(u)$ provided $\delta_\alpha,\delta_\beta={\cal O}(u)$, since only multiplications are involved in the backward relations.
\hfill \qed 

This lemma implies that relative small perturbations in the coefficients of $q(y)$ can be mapped to relative small perturbations in the coefficients of $p(x)$, both element-wise and norm-wise.

\subsection{Calculating the roots}

 The roots of the polynomial $q(y):=y^2-2\beta y +e$ are given by
$$ y_1=\beta + \sqrt{\beta^2-e}, \quad y_2=\beta - \sqrt{\beta^2-e}.
$$
 The way that these roots are computed depend now on the values of $\beta$ and $e$.
\begin{enumerate}
\item[Case 1:] $e=-1$ (real roots) 
$$y_1=fl(\beta + fl(\sqrt{fl(\beta^2+1)}), \quad y_2=-fl(1/y_1). $$
\item[Case 2:] $e=1$ and $\beta \ge 1$ (real roots)
$$y_1=fl(\beta + fl(\sqrt{fl(\beta+1)(\beta-1)}), \quad y_2=fl(1/y_1). $$
\item[Case 3:] $e=1$ and $\beta < 1$ (complex conjugate roots)
$$y_1=\beta + \j fl(\sqrt{fl(\beta+1)(1-\beta)}), \quad y_2=\overline{y}_1. $$
\end{enumerate}

Let us now check that the roots are computed in a forward stable manner. The error analysis
for the operations performed in the IEEE floating point standard give the following bounds.

\begin{enumerate}
\item[Case 1:] $e=-1$ (real roots) 
$$\hat y_1=(\beta + \sqrt{(\beta^2+1)})(1+\eta_3), \quad \hat y_2=-(1/\hat y_1)(1+\eta_1), \quad |\eta_i|\le \gamma_i. $$
\item[Case 2:] $e=1$ and $\beta \ge 1$ (real roots)
$$\hat y_1=(\beta + \sqrt{(\beta+1)(\beta-1)})(1+\eta_4), \quad \hat y_2=(1/\hat y_1)(1+\eta_1), \quad |\eta_i|\le \gamma_i. $$
\item[Case 3:] $e=1$ and $\beta < 1$ (complex conjugate roots)
$$\hat y_1=\beta + \j(\sqrt{(\beta+1)(1-\beta)})(1+\eta_3), \quad \hat y_2=\overline{\hat y}_1, \quad |\eta_i|\le \gamma_i. $$
\end{enumerate}

Notice that these bounds imply forward stability for all these computations. Combining this with Lemma 1, we have thus shown the following theorem.
\begin{theorem}
The computed roots $\hat y_i, i=1,2$ of the polynomial $q(y)$ satisfy the relative forward bounds
$$ |\hat y_1-y_1|\leq \Delta |y_1|, \quad |\hat y_2-y_2|\leq \Delta |y_2|, \quad \Delta={\cal O}(u),$$
and the transformed roots $\hat x_i=fl(-\alpha \hat y_i), i=1,2$ satisfy the mixed bounds
\begin{center}
$ a(x-\tilde x_1)(x-\tilde x_2)= ax^2+\hat b x + \hat c $ \\ \medskip
$ |\hat x_1-\tilde x_1| \leq \Delta |\hat x_1|, \;  |\hat x_2-\tilde x_2| \leq \Delta |\hat x_2|,$  \\ \medskip
$ |b-\hat b| \leq \Delta |b|, \;  |c-\hat c| \leq \Delta |c|, \; \Delta = {\mathcal O}(u).$
\end{center}
\end{theorem}

We can therefore also evaluate the backward bound by recomputing the sum and product of the computed roots. We first point out that the sum and product will be real because even when the two computed roots $\hat y_1$ and $\hat y_2$ are complex they will be exactly complex conjugate.

Since the product of the exact roots is $e=\pm 1$, and the computed roots are forward stable, we obviously have that the product of the computed roots satisfies
$$ \hat y_1 \hat y_2=e(1+{\cal O}(u))$$
which is element-wise backward stable in a relative sense.

For the sum of the computed roots, it is more problematic.  Since $|y_1|\ge |y_2|$ and both these roots are computed in a
forward stable way, we will have that
\begin{equation} \label{sum}  \hat y_1 + \hat y_2=\beta +{\cal O}(u)\hat y_1\end{equation}
but $\hat y_1$ can be much larger than $\beta$ and the backward error will then be much larger than $\beta \cdot {\cal O}(u)$.
Let us analyze the three cases. For Case 3 the sum of the computed roots is {\em exactly} $2\beta$ since this is a representable number. In Case 2, $\hat y_1\le 2\beta$ and \eqref{sum} then implies backward stability for the element $\beta$. But when $\beta\ll 1$ we can {\em not} obtain a sufficiently small backward error for \eqref{sum} since the recomputed sum has an error that is of the order of ${\cal O}(u)\hat y_1 \gg {\cal O}(u)\beta$. It is in this special case that element-wise backward stability gets lost.

\section{Complex coefficients}

The case where $b$ and/or $c$ are zero are again easy to handle but the relative error bounds
are slightly larger. Since exact error bounds are more difficult to describe, we preferred to just indicate their order of magnitude. 
 Let us first treat the case of 0 values.

If $c=0$, then the roots can be computed as follows
$$ x_1:=-b/a, \quad x_2=0$$
which is element-wise backward stable since under the IEEE floating point standard, we have that the {\em computed} roots satisfy (see \cite{Higham})
$$ \hat x_1= -fl(b/a)=-b(1+\delta)/a=-\hat b/a, \quad \hat x_2 = 0, \quad |\delta| = {\cal O}(u).$$
The backward error then indeed satisfies the relative element-wise bounds
$$ |\hat b - b| \leq |\delta| |b|, \quad |\hat c - c| \leq 0|c|, \quad |\delta| = {\cal O}(u). $$

If $b=0$ then the roots can be computed as follows
$$ x_1 = \sqrt{-c/a}, \quad x_2 := -x_1,$$
which is also element-wise backward stable since under the IEEE floating point standard, we have that the  {\em computed} roots satisfy (see \cite{Higham})
$$ \hat x_1=fl(\sqrt{fl(-c/a)})= \sqrt{-c(1+\eta)/a}, \quad \hat x_2=-\hat x_1, \quad |\eta| = {\cal O}(u).
$$
The backward error then satisfies the relative element-wise bounds
$$ |\hat b - b| \leq 0|b|, \quad |\hat c - c| \leq |\eta| . |c|, \quad |\eta| = {\cal O}(u). $$

\medskip

When there are no zero values, we again first apply a scaling of the problem.

\medskip

{\bf Scaling the polynomial $p(x)$}\\  
As in the real case, we scale the coefficients as follows~: $b_1:=b/a, c_1:=c/a$,
which can be performed in a backward and forward stable way since $a\neq 0$. 
According to the IEEE floating point standard we have indeed that
$$ |b-b_1|\leq |\Delta||b| , \; |c-c_1|\leq |\Delta||c|, \quad  |\Delta|={\cal O}(u). $$
This implies that we can as well look at the monic polynomial \\ $p(x)/a=p_1(x)=x^2+b_1x+c_1$.

\medskip

{\bf Scaling the variable $x$}\\ 
This becomes more complicated for the case of complex coefficients. We now have that $y:=-x/\alpha$ where 
$|\alpha|:=\sqrt{|c_1|}$ and $\arg(\alpha)=\arg(b_1)$.
This implies that we can consider again the polynomial 
\begin{equation}\label{y} q(y)=y^2 -2 \beta y +e =0, \end{equation}
where $\beta \in \RR_+$ and $|e|=1$. The formulas to compute $\alpha$, $\beta$ and $e$ are 
$$ b_1=|b_1|e_b, \; c_1=|c_1|e_c, \; \alpha:=e_b\sqrt{|c_1|}, \; \beta := |b_1|/(2\sqrt{|c_1|}), \; e:=e_c/(e_b)^2$$
where $e_b:=\arg(b_1)$ and $e_c:=\arg(c_1)$. For computational reasons, we will also compute the square root 
$f$ of $e$, i.e. $f^2=e$.

\medskip

We have again a similar lemma describing the transformation between the coefficients of the polynomials
\begin{lemma} The transformations 
$$ [\alpha,\beta,f]=h_a[b,c] \quad \mathrm{and} \quad [b,c]=h_a^{-1}[\alpha,\beta,f]$$ 
between the polynomial $p(x)=ax^2+bx+c, \; a\neq 0$ and the monic polynomial $p(-\alpha y)/(a\alpha^2)=q(y)=y^2-2\beta y +f^2$ defined by the forward and backward relations 
$$ \alpha:=\mathrm{arg}(b/a)\sqrt{|c/a|}, \quad \beta := |b/a|/(2\sqrt{|c/a|}), \quad f=\sqrt{arg(b/a)}/arg(c/a)
$$
and
$$ b=-2a\beta\alpha, \quad c=af^2\alpha^2,
$$
where $a$ is not perturbed, are both element-wise well-conditioned maps.
\end{lemma}
{\bf Proof}. 
The proof is very similar, except for the fact that the quantities are complex, except for $\beta$ which is real, and $f$ that can be parameterized by a real angle.
\hfill \qed 

This lemma implies again that relative small perturbations in the coefficients of $q(y)$ can be mapped to relative small perturbations in the coefficients of $p(x)$, both element-wise and norm-wise.

{\bf Calculating the roots}  \\
 The roots of the polynomial \eqref{y} are now given by
$$ y_1=\beta + \sqrt{\beta^2-f^2}, \quad y_2=\beta - \sqrt{\beta^2-f^2}.
$$
But we need only consider the case where $e=f^2$ is not real since otherwise we can apply the analysis of the previous section. The algorithm for computing the two roots is to first compute $y_1$ as the root of largest module, and then to compute $y_2$ using $y_2=f^2/y_1$. 
If we compute the square root of the complex number $\beta^2-f^2$  as $$\gamma=\sqrt{(\beta-f)(\beta+f)}$$
then the roots are given by
$$ y_{1}:= \beta +  \mathrm{sign}(\mathrm{real}(\gamma))\gamma, \quad y_2=f^2/y_1.
$$
The rounding errors can be written as follows
\begin{center}
$ \hat \gamma = \sqrt{\beta^2-f^2}(1+\delta_1)$\\ \medskip
$\hat y_1 = (\beta + |\mathrm{real}(\hat \gamma)|) (1+\delta_2) + \j  \mathrm{sign}(\mathrm{real}(\hat\gamma))\mathrm{imag}(\hat\gamma),$
\\ \medskip
$\hat y_2 = f^2(1+\delta_3)/\hat y_1,$
\end{center}
where all $|\delta_i|,i=1,2,3$ are of the order of the unit round-off $u$.
These formulas yield that $y_1$ and $y_2$ can be computed in a forward stable way.

The backward error analysis of these operations will be a problem when $\beta$ is much smaller than $|f|$.
This leads to the same conclusions as in the case of real coefficients: when the sum of the roots is much smaller than
the roots themselves, the relative backward error on the sum can be large, despite the fact that the forward errors 
on the computation as a function of $\beta$ and $f$ are small. 

\section{Comparing the different stabilities}

In this section we compare the different types of stability in terms of the constraints that they impose
on the computed roots.
First of all, it is obvious that EBS implies EMS since EMS follows from EBS by just choosing
$$ \tilde x_1 = \hat x_1, \quad \mathrm{and} \quad \tilde x_2 = \hat x_2. 
$$ 

We now prove that EBS implies NBS, which is slightly more involved.
\begin{lemma}
Let the computed roots $\hat x_1$ and $\hat x_2$ of $p(x)=ax^2+bx+c$ satisfy 
\begin{center} 
$ a(x-\tilde x_1)(x-\tilde x_2)= ax^2+\hat b x + \hat c $ \\ \medskip
$ |\hat x_1-\tilde x_1| \leq \Delta |\tilde x_1|, \;  |\hat x_2-\tilde x_2| \leq \Delta |\tilde x_2|,$  \\ \medskip
$ |b-\hat b| \leq \Delta |b|, \;  |c-\hat c| \leq \Delta |c|, \; \Delta = {\mathcal O}(u),$
\end{center}
then they also satisfy the norm-wise bound
\begin{center}
$ a(x-\hat x_1)(x-\hat x_2)= ax^2+\hat {\hat b} x + \hat {\hat c} $ \\ \medskip
$ \left\| \left[ \begin{array}{ccc} a & b & c  \end{array}\right] -\left[ \begin{array}{ccc}  a &  \hat {\hat b} & \hat {\hat c} \end{array}\right ]\right\| \leq  3 \Delta \left\| \left[ \begin{array}{ccc}  a & \hat b & \hat c \end{array}\right]\right\| , \quad \Delta = {\mathcal O}(u).$
\end{center}
\end{lemma}
{\bf Proof}\\
It follows from the EMS constraints that 
$$ \hat {\hat b} = \hat b + a(\hat x_1 - \tilde x_1 + \hat x_2 - \tilde x_2) , \quad \mathrm{and} 
\quad \hat {\hat c}=\hat c + a(\tilde x_1\tilde x_2 - \hat x_1\hat x_2),
$$
which yields the bounds
$$ |b- \hat {\hat b}| \leq |b-\hat b| + |a|(|\hat x_1 - \tilde x_1|+|\hat x_2 - \tilde x_2|)
 , \quad \mathrm{and} \quad
|c- \hat {\hat c}| \leq |c-\hat c| + |a|(|\tilde x_1\tilde x_2 -\hat x_1\hat x_2|).$$
Using the constraints of EMS we then also obtain
$$ |b- \hat {\hat b}| \leq \Delta|b| + \Delta|a|(|\tilde x_1|+|\tilde x_2|))
 , \quad \mathrm{and} \quad
|c- \hat {\hat c}| \leq \Delta|c| + \Delta|a|(|\tilde x_1| |\tilde x_2|).$$
Switching to norms and using the triangle inequality then yields
$$ \left\|\left[ \begin{array}{ccc} \; 0, \; & |b-\hat{\hat b}|, & |c-\hat{\hat c}| \end{array}\right]\right\|_2 \leq 
\Delta \left\|\left[ \begin{array}{ccc} \; 0,  & |{\hat b}|, & |{\hat c}| \end{array}\right]\right\|_2+ 
\Delta |a|\left\|\left[ \begin{array}{ccc} \; 0, \; & |\tilde x_1|+|\tilde x_2|, & \; |\tilde x_1\tilde x_2| \; \end{array} \right]\right\|_2. $$
Because of Lemma 4 in the appendix we also have 
$$ 
\Delta |a|\left\|\left[ \begin{array}{ccc} \; 0, \; & |\tilde x_1|+|\tilde x_2|, & \; |\tilde x_1\tilde x_2| \; \end{array} \right]\right\|_2 
\leq \sqrt{3}\Delta \left\|\left[ \begin{array}{ccc} \; a,  & |{\hat b}|, & |{\hat c}| \end{array}\right]\right\|_2 $$
and we finally obtain the norm-wise bound 
$$ \left\|\left[ \begin{array}{ccc} \; 0, \; & |b-\hat{\hat b}|, & |c-\hat{\hat c}| \end{array}\right]\right\|_2 \leq 
3 \Delta \left\|\left[ \begin{array}{ccc} \; a,  & {\hat b}, & {\hat c} \end{array}\right]\right\|_2. $$
\hfill \qed

We then need to show that in general, EBS can not always be satisfied, i.e. there does NOT exist any algorithm that achieves this.
A counterexample is given by the polynomial $y^2-2\beta -1$ where $\beta= 2^{-t}+2^{-2t}$ and $2^{-2t}\le u/2$ while  $2^{-t}\approx\sqrt{u}$.
One easily checks that $\beta$ is a representable number and that the roots of the polynomial are given by the expansion
$$ y_1= 1 + \beta + \beta^2/2 - \beta^4/8 + ..., \quad y_2= -1 + \beta - \beta^2/2 + \beta^4/8 + ...$$
Their exactly rounded values are given by the representable numbers
$$ \hat y_1 = 1 +  2^{-t}, \quad \hat y_2 = -1 +  2^{-t}
$$
which gives a sum equal to the representable number  
$$ \hat y_1 + \hat y_2 = 2. 2^{-t},
$$
but that yields a relative error of the order of $\sqrt{u}$ !  Moreover, all other representable numbers in the neighborhood of $y_1$ and $y_2$ 
are on a grid of size $u$ and all possible combinations of their sums will still have a comparable relative error. 
It is thus impossible to find representable numbers that would satisfy the EBS property.

\section{Numerical results}

\begin{center}
\includegraphics[scale=1]{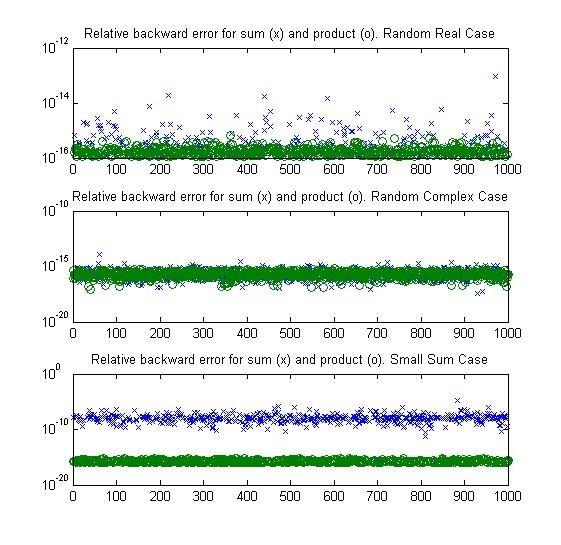}
\end{center}

We tested this routine for the relative backward errors on three sets of 1000 random quadratic polynomials.
We first took random real polynomials, then random complex polynomials, and finally random real polynomials with a very small 
sum (of the order of $\sqrt{\epsilon})$. The test results are given below.

The first plot clearly shows EBS, since the relative errors of the recomputed sums and products of the roots is of the order of 
the unit round-off $u$. The second plot shows the same results for polynomials with complex coefficients. The third plot shows that 
for real polynomials $q(y)$ with a very small (but non-zero) coefficient $\beta$, EBS can not be ensured by our algorithm. This is consistent with our analysis that shows that there does not exist any algorithm to ensure EBS for such polynomials.

\section*{Appendix A}

\begin{lemma}
For any real numbers $a$, $b$ and $c$ we have the inequality
$$ (|a|+|b|)^2 +|ab|^2 \leq 2c^2 +(a+b)^2 +(1+2/c^2)(ab)^2 $$
which also implies for $c^2=3/2$ that
$$ \left\|\left[ \begin{array}{ccc} \; 0, \; & |a|+|b|, & |ab| \end{array}\right]\right\|_2^2 \leq 3
\left\|\left[ \begin{array}{ccc} \; 1, \; & (a+b), & \; ab \; \end{array} \right]\right\|_2^2 . $$
\end{lemma}
{\bf Proof}\\
The first inequality follows from 
$$ (|a|+|b|)^2 = (|a|-|b|)^2 +4|ab| \leq (a+b)^2 + 4|ab|,
$$
and 
$$ 2(c-|ab|/c)^2\ge 0 \quad \Longrightarrow \quad 4|ab| \leq 2c^2 + 2(ab)^2/c^2. $$

\section*{Appendix B}

\begin{verbatim}
function [x1,x2,beta,e,scale] = quadroot(a,b,c)
% Function [x1,x2,beta,e] = quadroot(a,b,c)  computes the two roots 
% x1 and x2 of a quadratic polynomial ax^2+bx+c=0 in a stable manner
beta=[];e=[];scale=[];
% special cases of zero elements
if a==0, return, else b1=b/a;c1=c/a; end
if b==0, x1=sqrt(-c1);x2=-x1; return, end
if c==0, x1=-b1; x2=0; return, end
% generic case 
if isreal([b1,c1]),
    % with real coefficients
    c1abs=abs(c1);
    scale=sqrt(c1abs)*sign(b1);
    beta=b1/(2*scale);
    e=sign(c1);
    % computing the roots
    if e==-1, y1=beta+sqrt(beta^2+1);y2=-1/y1;
    else, 
        if beta >= 1, y1=beta+sqrt((beta+1)*(beta-1)); y2=1/y1;
        else, im=sqrt((beta+1)*(1-beta));y1=beta+j*im;y2=beta-j*im;
        end
    end
 else,
    % with complex coefficients
    scale=sign(b1)*(sqrt(abs(c1)));
    beta=abs(b1)/(2*sqrt(abs(c1)));f=sqrt(sign(c1))/sign(b1);
    gamma=sqrt((beta-f)*(beta+f));
    y1=beta+sign(real(gamma))*gamma;
    y2=f^2/y1;
end
x1=-y1*scale;x2=-y2*scale;
\end{verbatim}

\end{document}